\documentclass[10pt]{amsart}
\usepackage{amsmath,amssymb,amscd,enumerate}
\usepackage{amsmath,amssymb,amscd,enumerate,bbm}
\usepackage{graphicx}
\usepackage{color}
\usepackage{mathbbol}

\renewcommand{\MR}[1]{}

\newcommand{\titl}{DERESONATING A TATE PERIOD}

\title{{\titl}}
\author{V. GOLYSHEV}
\date{}

\newcommand{\cal}{\mathcal}

\def\R{{\Bbb{R}}}
\def\C{{\Bbb{C}}}
\def\Q{{\Bbb{Q}}}
\def\Z{{\Bbb{Z}}}

\def\HH{{\cal{H}}}
\def\DD{{\cal{D}}}

\def \DD {\mathcal{D}}

\renewcommand{\phi}{{\varphi}}

\newcommand{\rk}{\mathop{\mathrm{rk}}}

\newcommand{\Pic}{\mathop{\mathrm{Pic}}}

\newcommand{\Spec}{{\text{Spec }}}

\newcommand{\chern}{\mathop{\mathrm{ch}}}

\newcommand{\LCM}{\mathop{\mathrm{LCM}}}

\newcommand{\tensor}{\otimes}

\newcounter{pphcounter}[section]
\renewcommand{\thepphcounter}{\thesection.\arabic{pphcounter}}
\newcommand{\pph}{\bigskip \refstepcounter{pphcounter}
    \bf  \thepphcounter. \rm}

\makeatletter
\@addtoreset{equation}{pphcounter}
\makeatother

\newcommand{\coro}{\bf Corollary. \rm}

\newcommand{\propo}{\bf Proposition. \rm}
\newcommand{\theo}{\bf Theorem. \rm}

\renewcommand{\proof}{\bigskip \bf Proof. \rm}

\newcommand{\benum}{\begin{enumerate}}
\newcommand{\eenum}{\end{enumerate}}

\def\A1{{{\Bbb{A}}^1}}
\def\P1{{{\Bbb{P}}^1}}

\def\Gm{{\bf G_m}}

\def\SL2{{\mathrm SL2}}

\renewcommand{\P}{{\Bbb{P}}}

\newcommand{\lto}{\longrightarrow}

\newcommand{\comb}[2]{\left( {#1} \atop {#2} \right)}

\newcommand{\ww}{\wedge ^2}
\newcommand{\bgam}{\mathbb{\Gamma}}
\newcommand{\hbullet} {H^{\scriptscriptstyle ^\bullet}}

\begin{document}

\begin{center}
\maketitle

\bigskip


\bigskip

\end{center}

\bigskip

\parbox{340pt}{\small \bf Abstract. \rm We introduce a technique to
compute monodromy periods in certain families of algebraic varieties by
perturbing (`deresonating')
the fiberwise Betti to de Rham comparison off the motivic
setting. As an application, we find Apery periods of Grassmannians $G(2,N)$
and identify the Apery
numbers for the equations D3 of the Mukai threefolds with  certain  $L$--values.
We show that the argument of the
$L$--function is $3$ for the rational and
$2$ for the non--rational Mukai threefolds.}

\bigskip
\bigskip
\bigskip

\section{Apery limits for equations D3 and deresonating}

How can one read the topology or geometry of a Fano variety off its
Landau--Ginzburg model? Van Enckevort and van Straten were able  \cite{ES06},
given the LGs of certain four--dimensional Fanos $F$, to recover the Chern classes
of the restriction of the tangent bundle to the
Calabi--Yau anticanonical section  in $F$. L.~Katzarkov has suggested recently
to study the birational type of  $F$ by looking at its Landau--Ginzburg model.
We study this question for the class of Fano 3--folds considered
by Sh. Mukai, namely, the complete intersections in the Grassmannians of
simple algebraic groups.
We link the rational type of the generic variety in a family
with a certain period of Tate(--Artin) type, the so called Apery constant.

\pph \bf Apery limits. \rm We shall say that a linear homogeneous
recurrence $R$ with polynomial coefficients is a recurrence of the Apery type,
if there is a Dirichlet  character with  $L$--function
$L(s), \; \text{ an argument } s_0 \in \Z, \: s_0 > 1$ and two solutions of  $R$
$$u(n)=a_n, b_n \in \Q$$ such that
$$ \lim_{n \to + \infty} \dfrac{b_n}{a_n}= c \, L(s_0), \text{ где } c \in  \Q^\ast.$$
We shall say that the limit above is an Apery limit of the recurrence $R$.

\bigskip

The solution spaces of the recurrences of the arithmetic type
we consider are typically endowed with two
filtrations. One is the Dwork type filtration by the rate of growth of the
denominators. The other indicates how many of the leading terms of the
solution annihilate. One can specify what $a_n$ and $b_n$
are in the presence of the two filtrations and
 speak of the \emph{Apery constants} in
such cases. In our situation,  $a_n$ is characterized simply
as the unique integral solution normalized by
$a_0=1$, and $b_n$, as the unique solution with
$b_0=0,\; b_1=1.$
We refer the reader to \cite{ASZ08} for intriguing numerical findings
on the Apery constants for recurrences that arise from certain differential equations
of order $4$ and $5$.

In this language, our result
\ref{constants-of-recurrences-for-threefolds}
says  that the quantum
recurrences (sec. \ref{constructing-quantum-recurrences})
of the
Mukai threefolds
$V_{10},V_{12},V_{14},V_{16},V_{18}$ are recurrences of the
Apery type, the respective Apery constants being
$\dfrac{1}{10}\, \zeta(2)$, $\dfrac{1}{6} \, \zeta(3)$,$\dfrac{1}{7} \, \zeta(2)$,
$\dfrac{7}{32} \, \zeta(3)$, $\dfrac{1}{3}  \, L(\chi_3,3)$.  We see, in particular,
that the argument of the
$L$--function is $3$ for the rational and
$2$ for the non--rational Mukai threefolds.

For the rational varieties,
one can compute the Apery constant by the method of Beukers
that uses modularity of
the regularized quantum differential equation~(sec.~\ref{rational-cases}).
For the non-rational varieties (sec. \ref{non-rational-cases})
we introduce a new method, deresonating.

\pph \bf Deresonating. \rm
The formula
$$\log \Gamma (1-t) = \gamma t + \sum_{i=2}^{\infty} \dfrac{\zeta(i)}{i}t^i$$
suggests treating a Tate motive as if it were
a resonance limit of a non--motivic
entity in a wider hypergeometric realm.
Put $D=t \dfrac{d}{dt}$, and let $\DD= \C [t, t^{-1}, D]$ denote
the ring of differential operators on the torus.
A hypergeometric D--module
$$
\HH(\alpha_i,\beta_j) = \DD /\DD \left( \prod_{i=1}^n (D-\alpha_i)-
t \prod_{j=1}^n (D-\beta_j)\right)
$$
with rational indices
$\alpha_i,\beta_j$
is motivic, i. e. is
a subquotient of a variation of periods
in a pencil of varieties over $\Gm$ defined over
$\overline{\mathstrut \Q}$. The local system of its solutions is endowed with
both Betti and de Rham structures.

No rational structure can exist in the case of irrational exponents, and yet,
according to Dwork, a motivic quantity in (a pullback of) a hypergeometric family should be extended along the space of hypergeometric indices.
A \emph{gamma structure} \cite{GM09}
 on a hypergeometric D--module
manifests itself as a rational
structure in the case of rational exponents and gives rise to an extension
of the Betti to de Rham comparison along the space of exponents,
i.e. in the non-motivic direction (hypergeometric non-periods).

One might attempt to use l'Hospital's rule to extract
certain Tate--type periods out of expressions in gamma values as follows:

1) realize a Tate motive in a degenerate limiting fiber
in a  family of hypergeometric pure motives;

2) perturb the hypergeometric exponents to
a non-resonant set;

3) pass to the degenerate nonresonant non-period matrix; compute it;

4) let the perturbation parameters tend back to $0$.

We call the process of so perturbing a Tate type period to an expression
in gamma--values \emph{deresonating}.

\bigskip

We will change slightly the proposed setup,
deresonating the Apery constants,
which are frequently periods in families related to hypergeometric families.
The Apery constants are monodromy periods as opposed to fiberwise periods.
The role of the Apery constants in topology can be explained in short as follows.  Quantum topology (\cite[4.2.1]{Dubrovin98}, \cite{Iritani07}, \cite{KKP08})
has discerned in the Todd genus (the topological embodiment
of what seems to be an atomic
thing, the logarithm of the multiplicative group law) the couple of \emph{gamma
genera}:
$$\frac{t}{1-\exp(-t)}=\exp(\frac{1}{2}t)\Gamma(1+\frac{t}{2\pi i})\Gamma(1-\frac{t}{2\pi i}).$$
so that, in particular, Riemann--Roch--Hirzebruch reads
$$\chi (A,B) = \int \chern A^* \, \gamma(X)^* \chern B \,\gamma(X)$$
for a Calabi--Yau $X$.
The individual gamma halves are invisible to the classical
Riemann--Roch--Hirzebruch but are seen by its quantum counterpart,
the \emph{Dubrovin conjecture} on the monodromy of the regularized
quantum DE of a Fano variety. Each Lefschetz submodule
in the
cohomology of a  Fano variety corresponds to a resonant summand
of the local quantum DE at infinity.
Finally, the monodromy of a DE
is related to the recurrence
on the expansion coefficients of its solutions via the Apery limits
by a classical argument of Beukers.
Thus, in Galkin's recent formulation
{\tt http://www.mi.ras.ru/ $\tilde{}$~galkin/work/zetagrass.pdf,}
a system of Apery constants
corresponds naturally to the Lefschetz decomposition of
the cohomology of a Fano, which leads to a definition of the
\emph{Apery class}. A conjecture put forward by Galkin and Iritani relates the Apery
class to the gamma class for such Fano varieties as Grassmannians.

\bigskip
By deresonating
the simplest Apery constant in the simplest family of Grassmannians,
$\, G(2,N)$'s,
we find (Theorem \ref{grassm}) its value to be
$\dfrac{6}{N^2(N+1)}\, \zeta (2)$ .
Formula
\ref{sketch-proof-grassm}(\ref{item-formula})
shows how the perturbed Apery constant, as an expression in hypergeometric exponents,
is assembled from the matrix of the base change between
the Frobenius basis and the gamma basis (which is a finite expression in gamma
values) and the limit ratio which is an algebraic expression in the hypergeometric
exponents (the \emph{sine formula}).
The sine formula is the simplest particular case of a more general
Vandermonde determinant formula; it corresponds to the choice of the second wedge as the
polynomial functor. It would be interesting to compare this method with the
methods of  \cite{Brown06}, \cite{Cartier02},\cite{GM04} where applicable.

\begin{center}
---
\end{center}

\bigskip
The very first Landau--Ginzburg models had
been studied by Beukers and Peters \cite{BP84} and Beukers and Stienstra \cite{SB85}
long before they were introduced in the
context of mirror symmetry. Namely, they showed that the recurrence that Apery
had used to prove irrationality of $\zeta (3)$ (resp. to find a measure
of irrationality of $\zeta (2)$) translated into Picard--Fuchs equation
in a family of $K3$ surfaces (resp. elliptic curves).  We have identified \cite{Golyshev07}
these families with the Landau--Ginzburg models of the Fano threefold $V_{12}$
(resp. del Pezzo surface of degree $5$):

\pph \bf Apery's recurrence for  $\zeta(3)$. \rm
Apery proved irrationality of $\zeta(3)$
in 1979 by considering the recurrence
$$n^3u_n-(34n^3-51n^2+27n-5)u_{n-1}+(n-1)^3 u_{n-2}=0.$$

Denote by $a_n$ he solution of the recurrence with $a_0=1,\; a_1=5$ and by
$b_n$ --- the solution that satisfies $b_0=0,\; b_1=1.$
Then \cite{MP05}

\benum

\item $\displaystyle \mid \zeta(3) - \frac{6 b_n}{a_n} \mid =
\sum_{k=n+1}^\infty \frac{6}{k^3a_ka_{k-1}}= o(a_n^{-2})$;

\item All $a_n$'s are integral; the denominator of $b_n$ divides
$12\LCM (1,2,\dots, n)^3$;

\item $a_n=O(\alpha^n)$ where  $\alpha$ is
 the root of the characteristic polynomial
 $x^2-34x+1$ that is greater in absolute value;

\eenum

Put   $\dfrac{6 b_n}{a_n}=\dfrac{p_n}{q_n}$ with
coprime integral  $p_n, q_n$. Then it follows from
   $\text{LCM} (1,2, \dots, n)\le (1+\epsilon) e^n$ that
$\displaystyle \mid \zeta(3) - \frac{p_n}{q_n} \mid=
o(q_n)^{(-1+\delta)}$ with some
$\delta > 0 .$~\footnote{One
can choose $\delta = \frac{\log \alpha -3}{\log \alpha +3}$.}
The key assertion here is (ii), which follows from the fact
that the solutions  $a_n$ и $b_n$
are iterated binomial sums:
$$a_n= \sum_{k=0}^n \comb{n}{k}^2 \comb{n+k}{k}^2$$
$$b_n= \frac{1}{6} \sum_{k=0}^n \comb{n}{k}^2 \comb{n+k}{k}^2
\left(\sum_{m=1}^n\frac{1}{m^2}+
\sum_{m=1}^{k} \frac{(-1)^{m-1}}{2m^3\comb{n}{m}\comb{n+m}{m}}\right).$$

\bigskip

Put $D=t \dfrac{\partial}{\partial t}$.
Denote $A(t)= \sum a_nt^n, \;B(t)=\sum b_nt^n$.
Put $L= {D}^{3}-t\left (2\,D + 1\right )\left (17\,{D}^{2}+17\,D+5\right )+{t}^{
2}\left (D+1\right )^{3}
.$
Then $LA=0$ and $(D-1)LB=0$.


\pph \bf Theorem of Beukers and Peters \cite{BP84} \rm .
  Assume that \mbox{$t \ne 0,1,(\sqrt{2} \pm 1)^4, \infty$}. Then:

\benum
\item The surface $S_t: 1-(1-XY)Z-tXYZ(1-X)(1-Y)(1-Z)=0$
is birationally equivalent to a K3 surface;

\item The form
$$\omega_t =  \left. \frac{dX \wedge dZ}{XZ(1-t(1-X)(1-Y)(1-ZY)} \right\vert_{S_t}$$
is the unique holomorphic 2--form on  $X_T$;

\item $\rk \Pic X_t \ge 19$ (and is  $19$ for generic $t$);

\item  The periods $y$ of the form $\omega_t$
satisfy the differential equation $Ly=0$.
\eenum

\rm
\qed

\medskip
We denote by $\mathrm{i}$
the square root of $-1$.

\section{Constructing quantum DEs and recurrences}
\label{constructing-quantum-recurrences}

Let $X$ be a Picard rank one Fano threefold. Denote by $-K$ the anticanonical
class of $X$.
Consider a one--dimensional torus
$\Gm=\Spec \C[t,t^{-1}]$.

\pph The following is the standard procedure to obtain quantum differential
equations and recurrences (cf e.g. \cite{Golyshev07}).

\medskip

\bf Step 1. \rm Define a trilinear functional $\left< \alpha, \beta, \gamma \right>$
on $\hbullet (X)$ setting
$$\left< \alpha, \beta, \gamma \right>=\sum t^d \cdot [\text{ number of maps } \mathbb{P} ^1 \longrightarrow X $$
$$
\text{ of degree $d$ with respect to $-K$ such that } 0  \text{ maps into a representative of } \alpha,$$
$$1  \text{ maps into a representative of } \beta,$$
$$ \infty \text{ maps into a representative of } \gamma ]
$$
One has:
$$\left< \;, \;, \right>:
(\hbullet (X)  \tensor \C[t,t^{-1}])^{\tensor \, 3} \longrightarrow \C[t, t^{-1}].$$

\bf Step 2. \rm Extend the Poincare pairing $(\, , \,)$ to the trivial vector bundle
$\HH = \hbullet (X) \tensor  \C[t, t^{-1}]$ horizontally.

\bf Step 3. \rm Turn the trilinear form into a multiplication law:
$$( \alpha \cdot \beta , \gamma ) =
\left< \alpha, \beta, \gamma \right>.$$

\bf Step 4. \rm Introduce a connection (= a D--module structure) in $\HH:$
for
$h \in \hbullet (X)=
\hbullet(X) \tensor 1 \subset  \hbullet (X) \tensor  \C[t, t^{-1}]$
one sets
$$Dh=-K \cdot h$$
(here $h$ is understood to be $h \tensor 1$).

\bf Step 5. \rm ``Convolute the system''
into a single scalar equation using
$1 \tensor 1$ for the cyclic vector:
$$\sum b_{ij}t^iD^j (1\tensor 1)=0.$$

\bf Step 6. \rm Translate this into a recurrence in $u$, the
expansion coefficients of the solutions,

$$R^\mathrm{irreg}: \sum b_{ij} u(n-j)(n-i)^j=0$$

\bf Step 7. \rm Pass to the equation $R^\mathrm{reg}$ whose solution is
 $u(n) n!$:
$$\sum_i u(n-i) \sum_j b_{ij}(n-i)^j
\rightsquigarrow
\sum_i u(n-i) \sum_j b_{ij}(n-i)^j \dfrac{n!}{(n-i)!}.$$



\pph \label{Mukai threefolds} We define the Mukai threefolds to be those Fano threefolds
with Picard rank $1$ that are complete intersections in the Grassmannians
of simple Lie groups other than projective spaces.
 They were considered by Sh. Mukai in \cite{Mukai92}.

\medskip


\begin{tabular}{l|l}

$V_{10}$  & a section of
the Grassmannian $G(2,5)$ by a quadric and a codimension 2 plane
\\  \hline
$V_{12}$  & a section of the orthogonal Grassmannian $O(5,10)$
by a codimension
7 plane \\  \hline $V_{14}$  & a section of the
Grassmannian $G(2,6)$ by a codimension 5 plane \\  \hline
$V_{16}$  & a
section of the lagrangian Grassmannian $L(3,6)$ by a codimension 3
plane \\  \hline
$V_{18}$  & a section of $G_2/P$ by a codimension 2
plane

\end{tabular}

\pph \bf The corresponding differential operators \rm
 \cite{Golyshev07}, \cite{Przyjalkowski07a}.

\begin{tabular}{l|l}

                                     \\\hline  $V_{10}$
& ${D}^{3}-2\,t\left (1+2\,D\right )
\left (11\,{D}^{2}+11\,D+3\right )-4\,
{t}^{2}\left (D+1\right )\left (2\,D+3
\right )\left (1+2\,D\right )

$

                                      \\\hline $V_{12}$
&${D}^{3}-t\left (1+2\,D\right )\left (
17\,{D}^{2}+17\,D+5\right )+{t}^{2}
\left (D+1\right )^{3}

$

                                    \\\hline   $V_{14}$
&${D}^{3}-t\left (1+2\,D\right )\left (
13\,{D}^{2}+13\,D+4\right )-3\,{t}^{2}
\left (D+1\right )\left (3\,D+4\right
)\left (3\,D+2\right )

$

                                     \\\hline  $V_{16}$
&${D}^{3}-4\,t\left (1+2\,D\right )
\left (3\,{D}^{2}+3\,D+1\right )+16\,{
t}^{2}\left (D+1\right )^{3}

$

                                      \\\hline $V_{18}$
&${D}^{3}-3\,t\left (1+2\,D\right )
\left (3\,{D}^{2}+3\,D+1\right )-27\,{
t}^{2}\left (D+1\right )^{3}

$
\end{tabular}

\pph  \label{constants-of-recurrences-for-threefolds}\theo
The quantum recurrences of the Mukai threefolds
$V_{10},V_{12},V_{14},V_{16},V_{18}$ are recurrences of the
Apery type.
The respective Apery constants are:

\medskip

\begin{tabular}{l|l|l|l|l}
 $V_{10}:$ & $V_{12}:$ & $V_{14}:$ & $V_{16}: $ & $V_{18}:$   \\ \hline
\rule[-2em]{0em}{4em} $\dfrac{1}{10}\, \zeta(2)$  & $\dfrac{1}{6} \, \zeta(3)$  & $\dfrac{1}{7} \, \zeta(2)$
 & $\dfrac{7}{32} \, \zeta(3)$   & $\dfrac{1}{3}  \, L(\chi_3, 3)$

 \end{tabular}

\section{$V_{12},V_{16},V_{18}$: rational cases}
\label{rational-cases}

\pph \label{Beukers's theorem} \bf Theorem of Beukers.
\cite[1.2]{Beukers87} \rm
Let $F(\tau)=\sum_{n=1}^\infty c_n q^n,\: q=e^{2\pi i \tau}$ be a modular form
of weight
$4$ and conductor $N$ which is Atkin--Lehner odd, i.e. satisfies
$$F(-1/N\tau)=- (\tau \sqrt{N})^4 F(\tau).$$

Put
$$f(\tau)=\sum_{n=1}^\infty \dfrac{c_n}{n^3} q^n.$$

Denote $L(F,s)=\sum_{n=1}^\infty \dfrac{c_n}{n^s}.$
Put $h(\tau)=f(\tau)-L(F,3)$. Then
$$h (\tau)= -(\tau \sqrt{N})^2 h(-1/N\tau).$$

\qed

\pph \bf ``Eisenstein harmonics''. \rm
We will need finite linear combination
of ``elementary Eisenstein series''
$$E_{2,i}(Q) = -\frac{1}{24}\,i \, (1 -  24 \sum_{n=1}^{\infty} \sigma (n)
Q^{in})$$
and
$$E_{4,i}(Q) = \frac{1}{240}\, i^2 (1 +  240 \sum_{n=1}^{\infty}
\sigma_3 (n) Q^{in}).$$

\bigskip

An implication of \ref{Beukers's theorem} is the following proposition worked out by
Beukers in the case $N=6$.

\pph \propo For the equations that correspond to the cases $V_{2N}$ with
$N=6,8,9$ the following hold:

\benum

\item The function $\sum  a_n t^n = \Phi (q(t))$
is an Atkin--Lehner odd weight $2$ modular form of level $N$, i.e.
$\Phi (-1/N \tau)= -N\tau^2 \Phi(\tau)$;

\item
$$\dfrac{\sum b_n(t(q))^n}{\sum a_n(t(q))^n} = \sum \dfrac{c_i}{i^3} q^i,$$
with the coefficients $c_i$ coming from a weight $4$ modular form $F$,
as in Beukers's theorem;

\item
$$\Phi (-1/N \tau)(f(-1/N \tau) -L(F,3))=  \Phi(\tau)(f (\tau) -L(F,3))$$

\item
the solution
 $$\sum (b_n - L(F,3) a_n) t^n =
 \Phi (\tau) h (\tau) $$
extends analytically beyond the singularity
 $\mathrm{i}/\sqrt{N}$.

\eenum

\proof

\benum

\item This was established in  \cite{Golyshev07}.

\item Straightforward. The expressions of $F$ in terms of the ``Eisenstein harmonics''
and the shape of $L$--function is the table below.

\vbox{
\hspace{0truemm}\makebox[520pt][l]{
\begin{tabular}{c|c}
\rule[-2em]{0em}{4em} Variety   & $\Phi$ \\ \hline
\rule[-2em]{0em}{4em}
$V_{12}$ & $ 5E_{2,1}-E_{2,2}+E_{2,3}-5E_{2,6} $  \\ \hline
\rule[-2em]{0em}{4em}
$V_{16}$ &   $ 4E_{2,1}-2E_{2,2}+2E_{2,4}-4E_{2,8} $\\ \hline
\rule[-2em]{0em}{4em}
$V_{18}$ &  $ 3E_{2,1}-3E_{2,9}$
\end{tabular}
}}

\medskip

\vbox{
\hspace{-20truemm}\makebox[520pt][l]{
\begin{tabular}{c|c|c}
\rule[-2em]{0em}{4em} Variety  & $F$  & $L(s)$ \\ \hline
\rule[-2em]{0em}{4em} $V_{12}$ &
$E_4-7 E_{4,2}+7E_{4,3}-E_{4,6}$  &
$(1-7 \cdot 2^{2-s} +7 \cdot 3^{2-s} - 6^{2-s})\zeta(s) \zeta(s-3)$ \\ \hline
\rule[-2em]{0em}{4em} $V_{16}$ &   $E_4-21/4 E_{4,2}+21/4E_{4,4}-E_{4,8}$
& $(1-21/4 \cdot 2^{2-s} +21/4 \cdot 4^{2-s} - 8^{2-s})\zeta(s) \zeta(s-3)$
\\ \hline
\rule[-2em]{0em}{4em} $V_{18}$ & $\sum_{n=1}^{\infty} (\dfrac{n}{3})\sigma_3(n) q^n$ & $\prod_p (1-(\dfrac{p}{3})p^{3-s})^{-1}(1-(\dfrac{p}{3})p^{-s})^{-1}$

\end{tabular}
}}

\medskip

\item Follows from  \ref{Beukers's theorem}.

\item Follows from (iii).

\eenum

\pph \coro One has for the varieties $V_{12},V_{16},V_{18}$
$$\lim_{n \to + \infty} \dfrac{b_n}{a_n} = L(F,3).$$

\proof The assertions that $\lim_{n \to + \infty} \dfrac{b_n}{a_n} = x$ and that
the solution
$ \sum {b_n}t^n - x \sum {a_n}t^n$ extends beyond the radius of convergence
of  $\sum {a_n}t^n$ and $\sum {b_n}t^n$ are equivalent.
One has only to note that the singularity with the smaller absolute value
is uniformized by the point  $\tau=\mathrm{i}/\sqrt{N}$ in each of these cases.

\section{$V_{10},V_{14}$: the non--rational cases}
\label{non-rational-cases}

\pph \theo \label{grassm} Let  $N$ be an integer $\ge 5$. Then the Apery
constant of the Grassmannian $G(2,N)$ is $\dfrac{6}{N^2(N+1)}\, \zeta (2)$:
there are two solutions of the regularized
quantum recurrence for  $G(2,N)$
$$ a_{Nn}^{G(2,N)}, \: b_{Nn}^{G(2,N)},$$
such that $a_n \in \Z,\, a_0=1,\, b_0=0, \, b_N=1$ и
$$ \lim_{n \to + \infty} \dfrac{b_{Nn}^{G(2,N)}}{a_{Nn}^{G(2,N)}}= \dfrac{6}{N^2(N+1)}\, \zeta (2).$$

\pph \coro One has:

\benum
\item the Apery constant for  $V_{10}$ is $\dfrac{1}{10}\, \zeta (2)$;
\item the Apery constant for $V_{14}$ is $\dfrac{1}{7}\, \zeta (2)$;
\eenum

{\bf Proof of the corollary.} \rm
i) The quantum Lefschetz theorem (cf e.g. \cite{Gathmann03})
implies that for the variety
$V_{10}$ one has
$$a_{n}^{V_{10}}=\dfrac{5}{2} a_{5n}^{G(2,5)}\dfrac{(n!)^3 (2n)!}{5n!} ,\:
b_{n}^{V_{10}}=b_{5n}^{G(2,5)}\dfrac{(n!)^3 (2n)!}{5n!}.$$

ii) Quantum Lefschetz says that for
$V_{14}$ one has
$$a_{n}^{V_{14}}={6}a_{6n}^{G(2,6)}\dfrac{(n!)^6}{6n!},\:
b_{n}^{V_{14}}=b_{6n}^{G(2,6)}\dfrac{(n!)^6}{6n!}.$$
\qed

\bigskip

The  proof of \ref{grassm} uses
(1) the representation of the quantum differential equation for a
Grassmannian as a (pullback of)
the wedge power of a hypergeomteric,
(2) deresonating, (3) computation \cite{GM09} of the monodromy
of a hypergeometric with respect to the \emph{gamma structure}, and
(4)  Dubrovin's description of the monodromy of the tensor power
of a semisimple Frobenius manifold.

\pph \label{sketch-proof-grassm} \bf Proof of Theorem \ref{grassm}. \rm
\benum
\item The quantum differential operator for $G(2,N)$ is
the second wedge of the $N$--Kummer pullback of the differential operator $(D-1/2)^N+t$
(which is, up to a convention, the second wedge of the q.d.o for the projective
space $G(1,N)$.
The second wedge generates the ideal
in the ring of the differential operators  $\C [t,D]$,
that annihilates all $2 \times 2$ minors of the fundamental matrix
of the q.d.o of the projective space).
This is a theorem of Bertram--Ciocan--Fontanine--Kim--Sabbah, \cite{BCK05}, \cite{KS08}.

\item Deresonate the hypergeometric  D--module as follows.
Introduce the operator
 $$L_{dr}=(D-1/2-u)(D-1/2+u)(D-1/2-e)(D-1/2+e)(D-1/2)^{N-4}+t.$$ Let
$S_u$ (resp. $S_e$) be the solution whose expansion starts with   $t^{1/2+u}$
(resp.$t^{1/2+e}$), and let $S_{-u}$ (resp. $S_{-e}$) be the solution whose
 expansion satrts with $t^{1/2-u}$(resp. $t^{1/2-e}$): put
$${ \bgam(n)=\frac{(-1)^n}{
\Gamma(1/2+e+n)\Gamma(1/2-e+n)\Gamma(1/2+u+n)\Gamma(1/2-u+n)\Gamma(1/2+n)^{N-4}}};$$
then
$$S_{-e}=\sum_{n=0}^\infty \bgam (1/2-e+n)t^{1/2-e+n}$$
$$S_{e}=\sum_{n=0}^\infty \bgam (1/2+e+n)t^{1/2+e+n}$$
$$S_{-u}=\sum_{n=0}^\infty \bgam (1/2-u+n)t^{1/2-u+n}$$
$$S_{u}=\sum_{n=0}^\infty \bgam (1/2+u+n)t^{1/2+u+n}$$

Let
$$R_u=S_{u}S_{-u}' - S_{u}'S_{-u} \text{ and } R_e= S_{e}S'_{-e} - S_{e}'S_{-e}$$
be the two solutions of $L_{dr}^{\ww}R=0$, and put
$R_u= \sum r_u^{(n)}t^n$, $R_e= \sum r_e^{(n)}t^n$. Then,
as we shall see in \ref{proof-formula-L}, the \emph{sine formula} holds
$$\lim_{n \rightarrow + \infty} \frac{r_e^{(n)}}{r_u^{(n)}}=
\frac{\sin (2 \pi  e)} {\sin (2 \pi  u)}  $$

\item \label{uppe} Put
$$A^{dr}(t)= \frac{R_e}{r_e^{(0)}},$$
$$B^{dr}(t)=
\frac{ r_e^{(0)}R_u-r_u^{(0)}R_e}
{r_e^{(0)}
r_u^{(1)}-r_u^{(0)}r_e^{(1)}}.$$
Then $A^{dr}(t)$ is a deresonation of
 $A(t)=\sum a_{Nn}^{G(2,N)} \dfrac{t^{n}}{(Nn)!}$,
and
$B^{dr}(t)$ is a deresonation of $B(t)={N} \sum b_{Nn}^{G(2,N)} \dfrac{t^{n}}{(Nn)!}$,
and the limit of the ratio of $n$--th respective coefficients is a perturbation
of the Apery constant of the Grassmannian:

$$\lim_{n \rightarrow + \infty}
\frac{ r_e^{(0)}r_u^{(n)}-r_u^{(0)}r_e^{(n)}} {r_e^{(0)}
r_u^{(1)}-r_u^{(0)}r_e^{(1)}}/
\dfrac{r_e^{(n)}}{r_e^{(0)}}
\rightarrow
{N} \lim_{n \rightarrow + \infty} \dfrac{b_{Nn}^{G(2,N)}}{a_{Nn}^{G(2,N)}}
\text{  as } u,e\rightarrow 0
$$

\item \label{item-formula} Combining it with the sine formula, we arrive at the following
expression for the perturbed Apery constant:
$$\text{p. A. c. }=\dfrac{1}{N} \cdot \frac{\sin(2\pi u)/\sin(2\pi e)\,r_e^{(0)2}-r_u^{(0)}r_e^{(0)}}
{r_u^{(1)}r_e^{(0)}-r_e^{(1)}r_u^{(0)}}
$$
It is a routine check that the limit of the p. A. c. as $u,e \to 0$ is  $\dfrac{\pi^2}{N^2(N+1)}$.

It only remains to check the sine formula.

\eenum

\pph \label{proof-formula-L} \bf Proof of the sine formula. \rm
It looks probable that the sine formula holds for any nonnegative $N$,
not necessarily integral.
\footnote{Don Zagier has suggested an appoach to the sine formula
which is based on the Poisson summation and does not require integrality of $N$.}
Let us give a sketch of a proof for even $N$.
We will use Dubrovin's extension of the Thom--Sebastiani
formula \cite{Dubrovin99}, \cite{Dubrovin04},
that expressses the monodromy of the so called second structural
connection of a product of two Frobenius manifolds in terms of the monodromy
data of the second structural connections of the factors.
This approach will work in a much wider framework: one can apply arbitrary
polynomial functors to Kloosterman type objects and compute the Apery limits
for the resulting DE's.

Identify  $\C$ with the real space
$\R ^2$, and fix a linear form $h: \R ^2 \lto \R$ in general position.
Consider the following objects of linear algebra: a vector space endowed with
a non--symmetric bilinear form $[\, , \,)$
and a choice of a
semiorthogonal basis
$\left< v_1, \dots, v_n \right>$ of $V$ compatible
with a marking $\mu = \mu_V: [1,\dots,n] \lto \C $:
\begin{eqnarray*}
\left[v_i,v_j \right) & = & 0
\text{ when } h(\mu(i)) > h(\mu(j)), \\
\left[ v_i,v_i \right) & = & 1
\text{ when }  1 \le i \le n.
\end{eqnarray*}
These may be used to produce polarized local systems on
 $\C \setminus \{ \mu(i) \}$
by identifying the fiber with $V$ ,
endowing it with the form $[\, , \,]$ (resp. $(\, , \,)$) ---
the (skew)symmetrization of the form $[,)$, choosing infinity for the base
point, joining it with the points $ \mu_i=\mu(i)$ with the level rays $h(x)=h(\mu_i)$
as paths and requiring that the turn around  $\mu_i$ act in the monodromy representation by the reflection w.r. to $v_i$. Vice versa, if there is such a local
system and the paths are given, one can consider the fiber with the basis of vanishing cycles and pass from the (skew)symmetrized form to the non--symmetric one
according to the order given by the values of $h(\mu_i)$.

These objects may be ``tensorized'': consider
$$V \tensor W, [\, , \, ), \left< v_1 \otimes w_1, \dots, v_n \otimes w_n \right>,  \mu$$
so that
$$[v_i \otimes w_j,v_k \otimes w_l]=[v_i,v_k)[w_j,w_l)$$
and
$$\mu (i,j)_{V \tensor W} = \mu_V(i)+\mu_W(j),$$
then passing to the (skew)symmetrized local system.

For simplicity, perturb our hypergeometric differential operator
$L^{dr}$
further to a non--resonance one
$$\prod_{i=1}^N (D-\alpha_i)+t;$$
assume also that every index  $\alpha_i$ occurs along with   $1-\alpha_i$.

\benum

\item
Put $t=-w^{-N}$, so that the coordinate $w$ is the Kummer pullback of
the coordinate $t$.
Clearly, wedging commutes with the Kummer pullback. We use the minus
sign for simplicity so as to deal with the roots of unity and not $-1$.
The final result is not affected, nor is it affected by the convention
to expand the local solutions around infinity and not zero.

\item The monodromy of the regularized
differential operator
$$ \prod_{i=1}^N (D+i)- w^N\prod_{i=1}^N (D+N\alpha_i)$$
can be described according to \cite{GM09} (cf also \cite[1.2]{Golyshev01})
as follows.
Put
$$H(y)= \frac{1-y^N}{\prod (1-y\exp (2\pi \mathrm{i} \alpha_i))},$$
and expand
$H(y)=1+\sum_{i=1}^\infty c_iy^i.$
Consider the $N$--dimensional $\C$--vector space  $\mathcal{V}$ with the basis $v_i$, endowed
with the symmetric bilinear form given by
\begin{eqnarray*}
\left( v_i,v_i \right) & = & 2
\text{ при }  1 \le i \le N, \\
\left(v_i,v_j\right) & = &  c_{\mid j-i \mid}
\text{ при } i \ne j.
\end{eqnarray*}

Set $U=\P ^1 \setminus \{ N \text{--th roots of unity }, \infty  \}$.
Interpret the monodromy representation as acting in the space
$\mathcal{V}, (\, , \,)$
so that the turn around  $\exp (2 \pi \mathrm{i} (j-1)/N)$ acts by the reflection
with respect to $v_j$.
Linear algebra shows that the expansion of the eigenvectors of the local monodromy
around $\infty$ with respect to the dual basis $\hat v_k$ is given by the Vandermonde
matrix (\cite{GM09}):
$$\mathbf{e}_i= \sum_j \exp( 2 \pi \mathrm{i} \alpha_i)^j \hat v_j.$$
Choose the eigenvectors
$\mathbf{e}_1, \mathbf{e}_2, \mathbf{e}_3, \mathbf{e}_4$
to correspond to the solutions
$S_{-e},S_{e},S_{-u},S_u$.

\item According to Dubrovin \cite[ch. 4, 5]{Dubrovin99},
the monodromy of the regularized second wedge of the operator
$$1-w^N\prod_{i=1}^N (D+N\alpha_i)$$
can be described as follows.
Let $V_{ij}=v_i \wedge v_j$ be the elements of the basis of the $\frac{N(N-1)}{2}$--
dimensional
$\C$--vector space $\Lambda^2 {V}$. Introduce the semiorthonormal form
$$[V_{ij},V_{kl})= [v_i,v_k)[v_j,v_l)-[v_i,v_l)[v_j,v_k),$$
and denote by $[\, , \, ]$ its skewsymmetrization.
The set of singularities of our DE is
$$U=\mathbb{C} \setminus \{\text{all sums of  pairs of distinct roots of unity}  \},$$
and the turn around $\exp (2 \pi \mathrm{i} (k-1)/N)+\exp (2 \pi \mathrm{i} (l-1)/N)$
acts by the reflection with respect to $V_{kl}$.
This assertion requires specifying the loops explicitly. In our situation,
a loop around
\mbox{$\mu_{12}=
\exp (2 \pi \mathrm{i} \cdot 0/N)+\exp (2 \pi \mathrm{i} \cdot 1/N)$}
can be chosen to be of the simplest shape: a ray from infinity to
$\mu_{12}$ --- a turn along a small circle around
$\mu_{12}$ --- the way back along the same ray.

\begin{figure}
  \includegraphics[scale=0.5,hiresbb=true,draft=false]{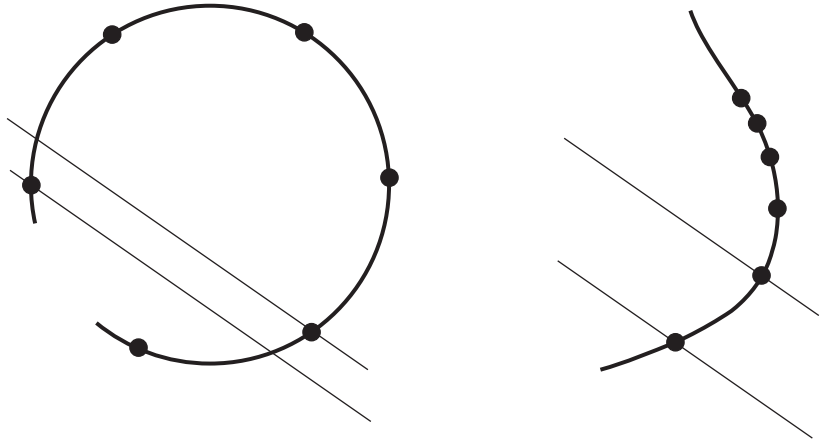}\\
  \caption{\small It is clear how to use the monodromy formula by shifting a little
  the singular points.
   The parallel lines are the level lines of $h$.}\label{fig1}
\end{figure}

The eigenvectors of the monodromy around $\infty$ have the
form
\cite[ch. 5]{Dubrovin99}
$$\mathbf{E}_{jk}= \mathbf{e}_j \wedge \mathbf{e}_k.$$
It is also clear that the solutions
$R^e,R^u$ correspond to
$\mathbf{E}_{12}, \mathbf{E}_{34}$
in the adopted notation.
Therefore the expansion coefficient at $\hat V_{12}$ in the expansion of
$\mathbf{E}_{12}$ with respect to the dual basis $\hat V_{ij}$
equals
$\sin (2 \pi  e)$,  and
the coefficient at $\hat V_{12}$ in the expansion of
$\mathbf{E}_{34}$ in $\hat V_{ij}$'s
equals
$\sin (2 \pi  u)$.

\item One notes now that the sine formula
$$\lim_{n \rightarrow + \infty} \frac{r_e^{(n)}}{r_u^{(n)}}=
\frac{\sin (2 \pi  e)} {\sin (2 \pi  u)}  $$
means that the solution ${\sin (2 \pi  u)} R^e- {\sin (2 \pi  e)} R^u$
can be extended analytically beyond the radius of convergence of
$R^e$ and $R^u$. This is equivalent to the coefficient at
$\hat V_{12}$
in the expansion of
${\sin (2 \pi  u)} R^e- {\sin (2 \pi  e)} R^u$ being
$0$, which has been proven above. This finishes the proof of the sine formula
and Theorem \ref{grassm}.

\eenum
\qed

\bigskip
\bigskip

I thank Sergei Galkin, Constantin Shramov and Duco van Straten
 for their remarks and suggestions.

\bigskip
\bigskip

\nocite{Andre04}
\nocite{KR07}


\begin{thebibliography}{BCFK05}

\bibitem[And04]{Andre04}
Yves Andr\'e, \emph{{An introduction to motives. Pure motives, mixed motives,
  periods. (Une introduction aux motifs. Motifs purs, motifs mixtes,
  p\'eriodes.)}}, {Panoramas et Synth\`eses 17. Paris: Soci\'et\'e
  Math\'ematique de France. xi, 261~p. EUR~26.00; \$~37.00 }, 2004 (French).

\bibitem[AvSZ08]{ASZ08}
Gert Almkvist, Duco van Straten, and Wadim Zudilin, \emph{{Ap\'ery limits of
  differential equations of order 4 and 5.}}, {Yui, Noriko (ed.) et al.,
  Modular forms and string duality. Proceedings of a workshop, Banff, Canada,
  June 3--8, 2006. Providence, RI: American Mathematical Society (AMS);
  Toronto: The Fields Institute for Research in Mathematical Sciences. Fields
  Institute Communications 54, 105-123 (2008).}, 2008.

\bibitem[BCFK05]{BCK05}
Aaron Bertram, Ionu\c t Ciocan-Fontanine, and Bumsig Kim, \emph{{Two proofs of
  a conjecture of Hori and Vafa.}}, Duke Math. J. \textbf{126} (2005), no.~1,
  101--136 (English).

\bibitem[Beu87]{Beukers87}
F.~Beukers, \emph{{Irrationality proofs using modular forms.}}, {Journ\'ees
  arithm\'etiques, Besan\c{c}on/France 1985, Ast\'erisque 147/148, 271-283
  (1987).}, 1987.

\bibitem[BP84]{BP84}
F.~Beukers and C.A.M. Peters, \emph{{A family of K3 surfaces and $\zeta(3)$.}},
  J. Reine Angew. Math. \textbf{351} (1984), 42--54 (English).

\bibitem[Bro06]{Brown06}
Francis~C.S. Brown, \emph{{Multiple zeta values and periods of moduli spaces
  $\bar {\mathfrak M}_{0,n}$. (P\'eriodes des espaces des modules $\bar
  {\mathfrak M}_{0,n}$ et valeurs z\^etas multiples.)}}, 2006, pp.~949--954.

\bibitem[Car02]{Cartier02}
Pierre Cartier, \emph{{Polylogarithmic functions, polyzeta numbers and
  pro-unipotent groups. (Fonctions polylogarithmes, nombres polyz\^etas et
  groupes pro-unipotents.)}}, {Bourbaki seminar. Volume 2000/2001. Expos\'es
  880-893. Paris: Soci\'et\'e Math\'ematique de France. Ast\'erisque 282,
  137-173, Exp. No. 885 (2002).}, 2002.

\bibitem[Dub98]{Dubrovin98}
Boris Dubrovin, \emph{{Geometry and analytic theory of Frobenius manifolds.}}
  (English).

\bibitem[Dub99]{Dubrovin99}
\bysame, \emph{{Painlev\'e transcendents in two-dimensional topological field
  theory.}}, {Conte, Robert, The Painlev\'e property. One century later. New
  York, NY: Springer. CRM Series in Mathematical Physics. 287-412 (1999).},
  1999.

\bibitem[Dub04]{Dubrovin04}
\bysame, \emph{{On almost duality for Frobenius manifolds.}}, {Buchstaber, V.
  M. (ed.) et al., Geometry, topology, and mathematical physics. Selected
  papers from S. P. Novikov's seminar held in Moscow, Russia, 2002--2003.
  Providence, RI: American Mathematical Society (AMS). Translations. Series 2.
  American Mathematical Society. 212. Advances in the Mathematical Sciences 55,
  75-132 (2004).}, 2004.

\bibitem[Gat03]{Gathmann03}
Andreas Gathmann, \emph{Relative {G}romov-{W}itten invariants and the mirror
  formula}, Math. Ann. \textbf{325} (2003), no.~2, 393--412. \MR{MR1962055
  (2004b:14096)}

\bibitem[GM04]{GM04}
A.B. Goncharov and Yu.I. Manin, \emph{{Multiple $\zeta$-motives and moduli
  spaces $\overline{\cal M}_{0,n}$.}}, 2004, pp.~1--14.

\bibitem[GM09]{GM09}
V.~Golyshev and A.~Mellit, \emph{Gamma structures and {G}auss's contiguity},
  2009.

\bibitem[Gol01]{Golyshev01}
V.~V. Golyshev, \emph{{Riemann-Roch variations.}}, Izv. Math. \textbf{65}
  (2001), no.~5, 853--881 (English. Russian original).

\bibitem[Gol07]{Golyshev07}
Vasily~V. Golyshev, \emph{{Classification problems and mirror duality.}},
  {Young, Nicholas (ed.), Surveys in geometry and number theory. Reports on
  contemporary Russian mathematics. Cambridge: Cambridge University Press.
  London Mathematical Society Lecture Note Series 338, 88-121 (2007).}, 2007.

\bibitem[Iri07]{Iritani07}
Hiroshi Iritani, \emph{Real and integral structures in quantum cohomology i:
  toric orbifolds}, 2007.

\bibitem[KKP08]{KKP08}
L.~Katzarkov, M.~Kontsevich, and T.~Pantev, \emph{Hodge theoretic aspects of
  mirror symmetry}, 2008.

\bibitem[KR07]{KR07}
C.~Krattenthaler and T.~Rivoal, \emph{{Hypergeometry and Riemann's zeta
  function. (Hyperg\'eom\'etrie et function z\^eta de Riemann.)}}, Mem. Am.
  Math. Soc. \textbf{875} (2007), 87 p. (English).

\bibitem[KS08]{KS08}
Bumsig Kim and Claude Sabbah, \emph{{Quantum cohomology of the Grassmannian and
  alternate Thom-Sebastiani.}}, Compos. Math. \textbf{144} (2008), no.~1,
  221--246 (English).

\bibitem[MP05]{MP05}
Yu.~I. Manin and A.~A. Panchishkin, \emph{{Introduction to modern number
  theory. Fundamental problems, ideas and theories. Transl. from the Russian.
  2nd revised ed.}}, {Encyclopaedia of Mathematical Sciences 49. Number Theory
  1. Berlin: Springer. xv, 514~p. }, 2005 (English).

\bibitem[Muk92]{Mukai92}
Shigeru Mukai, \emph{{Fano 3-folds.}}, {Complex projective geometry, Sel. Pap.
  Conf. Proj. Var., Trieste/Italy 1989, and Vector Bundles and Special Proj.
  Embeddings, Bergen/Norway 1989, Lond. Math. Soc. Lect. Note Ser. 179, 255-263
  (1992).}, 1992.

\bibitem[Prz07]{Przyjalkowski07a}
V.~V. Przyjalkowski, \emph{{Gromov-Witten invariants of Fano threefolds of
  genera 6 and 8.}}, 2007, pp.~433--446.

\bibitem[SB85]{SB85}
Jan Stienstra and Frits Beukers, \emph{On the {P}icard-{F}uchs equation and the
  formal {B}rauer group of certain elliptic {$K3$}-surfaces}, Math. Ann.
  \textbf{271} (1985), no.~2, 269--304. \MR{MR783555 (86j:14045)}

\bibitem[vEvS06]{ES06}
Christian van Enckevort and Duco van Straten, \emph{{Monodromy calculations of
  fourth order equations of Calabi-Yau type.}}, {Yui, Noriko (ed.) et al.,
  Mirror symmetry V. Proceedings of the BIRS workshop on Calabi-Yau varieties
  and mirror symmetry, December 6--11, 2003. Providence, RI: American
  Mathematical Society (AMS); Somerville, MA: International Press. AMS/IP
  Studies in Advanced Mathematics 38, 539-559 (2006).}, 2006.

\end{thebibliography}

\end{document}